\documentclass[12pt]{article}

\usepackage{amsmath}
\usepackage{amssymb}
\usepackage{amscd}
\usepackage{epsfig}

 \newcommand{\T}{\mathbb T}
 
 \newcommand{\X}{\mathbb X}
  \newcommand{\A}{\mathbb A}
 \newcommand{\qed}{\hfill $\square$}
 \newcommand{\bee}{\begin{equation}}
 \newcommand{\eee}{\end{equation}}
  
 \newcommand{\Lb}{\mbox {\boldmath ${\Lambda}$}}  
 \newcommand{\Gb}{\mbox {\boldmath ${\Gamma}$}}

 \newcommand{\Lbt}{\mbox{\tiny\boldmath ${\Lambda}$}}
 \newcommand{\Gbt}{\mbox{\tiny\boldmath ${\Gamma}$}}
 \newcommand{\Pb}{\mbox {\bf P}}

 \newcommand{\supreme}{\mbox{\rm sup}}

\newcommand{\be}{\begin{eqnarray}}
\newcommand{\ee}{\end{eqnarray}}
\newcommand{\supp}{\mbox{\rm supp}}

\newcommand{\dens}{\mbox{\rm dens}}

\newcommand{\R}{{\mathbb R}}

\newcommand{\Z}{{\mathbb Z}}

\newcommand{\Nat}{{\mathbb N}}

\newcommand{\Ak}{{\mathcal A}}

\newcommand{\Lam}{{\Lambda}}

\newcommand{\Gam}{\Gamma}

 \newtheorem{theorem}{Theorem}[section]
 
 \newtheorem{prop}[theorem]{Proposition}
 \newtheorem{cor}[theorem]{Corollary}

 \newtheorem{defi}[theorem]{Definition}

\numberwithin{equation}{section}

\begin{document}

\bibliographystyle{unset}

\title{{\sc A Characterization of Model Multi-colour Sets}}
\date \today 
\maketitle

 \centerline {
{\sc Jeong-Yup  Lee } and 
{\sc Robert V.\ Moody} 
}

{\small
\hspace*{4em}
  
\smallskip
\hspace*{4em} 
Department of Mathematics and Statistics, University of Victoria, \\
\hspace*{4em}
\hspace*{2.5em} Victoria, BC, V8W 3P4, Canada. 
\vspace{7mm}

\abstract Model sets are always Meyer sets, but not vice-versa. This article is about characterizing model sets (general and regular) amongst the Meyer sets in terms of two associated dynamical systems. These two dynamical systems describe two very different topologies on point sets, one local and one global. In model sets these two are strongly interconnected and this connection is essentially definitive. The paper is set in the context of multi-colour sets, that is to say, point sets in which points come in a finite number of colours, that are loosely coupled together by finite local complexity. 
\vspace{15 mm}

\section{Introduction}

Model sets, or cut and project sets, are an essential ingredient in the study of long-range aperiodic order and in both the theoretical and experimental sides of quasicrystal theory. Model sets are Delone subsets of space and under mild assumptions display such nice properties as finite local complexity, repetitivity, and pure point diffractivity while at the same time being free of any translational symmetry. However, it is not easy to characterize these sets. By definition a model set, say in real space $\R^n$, is the projection of part of a lattice from some
`super-space' $\R^n \times H$, where $H$ is some locally compact Abelian group, the part to be projected being determined by some suitable compact set $W$ in $H$. For a model set given to us in $\R^n$, all this additional baggage is hidden. The process of reconstructing it is part of the problem.

In the study of aperiodic order one needs to decide what it means for two aperiodic sets to look like each other. There are two very distinct ways in which this is commonly done. The first is to compare the two sets locally and see if, perhaps after a small translational shift, they agree on some large patch of space, say on a large ball around the origin. This leads to a topology on discrete sets which we call here the {\em local topology}. On the other hand one may take the view that it is not the local structure that is important, but rather the average structure. This leads us to determine that two discrete sets are close if, again perhaps after a small translational shift, the density of the symmetric difference of the two sets is small. This leads to the {\em autocorrelation topology}. 

Starting with a discrete set $\Lb$ one forms the translational orbit $\R^n + \Lb$.
The closures of this orbit in the local and autocorrelation topologies lead to two dynamical hulls,
$\X(\Lb)$ and $\A(\Lb)$. There is no reason to expect there to be any relationship between these two hulls since there is no inherent relationship between the local and autocorrelation topologies. But for regular model sets there is. In fact there is a continuous surjective $\R^n$-map $\beta:\X(\Lb) \longrightarrow \A(\Lb)$, the so-called {\em torus parameterization} \cite{martin}. This remarkable confluence of the two topologies seems to be a fundamental property of quasi-crystal theory. In this paper we show that, in the context of Meyer sets (see below), it characterizes the model sets.

The setting of the paper is not just sets, but more generally multi-colour sets, or multi-sets as we will refer to them henceforth.\footnote{Lagarias and Wang \cite{lawa} use the term multi-set in a different way, namely for point measures with non-negative integral coefficients.} These are sets in which points are coloured by some finite number $m$ of colours, so that in effect we have $m$ different types of points. This generalization is rather natural for any attempts to model physical structures by using points to designate atomic positions, and it is likewise essential in the study of aperiodic tilings, in which we represent tiles of different types by suitably coloured points. It is particularly relevant for substitution tilings where the substitution process can only be formulated at the level of point sets by the use of multi-sets.

The second generalization we make is to replace $\R^n$ by any compactly generated locally compact Abelian group $G$. This produces virtually no additional complications. 

Briefly the main theorem of the paper can be stated in the form:

\begin{theorem} 
Let $\Lb$ be a repetitive Meyer multi-set in $G$. Then there is a continuous $G$-map $\beta : \mathbb{X}(\Lb) \rightarrow \mathbb{A}(\Lb)$ which 
is one-to-one a.e. with respect to $\mathbb{A}(\Lb)$ if and only if $\Lb$ (or equivalently each element of $\mathbb{X}(\Lb)$) is a regular model multi-set.
\end{theorem}

This paper derives much from a longer work \cite{BLM}. That work is primarily
concerned with characterizing model sets in terms of dynamical
systems and their spectra, which is a considerably more complex story. Here our emphasis is to 
give a relatively short self-contained characterization of model sets in terms 
of the `coincidence' of two very different topologies, and to do this not
only for model sets, but also for model multi-sets.

The present paper also relies on the earlier work of Martin
Schlottmann \cite{martin}, and in fact the last section is essentially 
a reworking of his results, although the hypotheses are slightly different and again 
we are generalizing to multi-sets. 
In dealing with model sets, the devil is truly in the details. Model sets 
are in principle easy to understand, but the virtually unlimited scope of the
defining window $W$ produces all sorts of complications, and there are various
definitions of model sets in the literature, depending on exactly which
conditions are imposed on $W$. One of our contributions here is to
find a slightly different (and basical more general) definition for 
model sets (Def.\,\ref{modelSetDef}), which provides
for a smoother characterization than is otherwise possible and allows 
an efficient and fairly self-contained proof of the characterization of 
model sets.
The underlying context of the characterization is that of Meyer sets
(Def.\,\ref{def-Meyer}).

\section{Preliminaries}
\noindent
Let $G$ be a compactly generated locally compact Abelian group. 
A {\em multi-set} or {\em $m$-multi-set} in $G$ is a 
subset $\Lb = \Lam_1 \times \dots \times \Lam_m 
\subset G \times \dots \times G$ \; ($m$ copies)
where $\Lam_i \subset G$. 
We also write 
$\Lb = (\Lam_1, \dots, \Lam_m) = (\Lam_i)_{i\le m}$.
Although $\Lb$ is a product of sets, it is convenient to think
of it as a set whose points come in various `types' or `colours', $i$ being the
colour of the points in $\Lambda_i$. 
A {\em cluster} of $\Lb$ is, by definition,
a family $\Pb = (P_i)_{i\le m}$ where $P_i \subset \Lambda_i$ is 
finite for all $i\le m$.
Many of the clusters that we consider have the form
$A\cap \Lb := (A\cap \Lambda_i)_{i\le m}$, for a compact set
 $A \subset G$.
There is a natural
translation $G$-action on the set of multi-sets and their clusters
in $\Lb$. 
The translate of a cluster $\mbox{\bf P}$ by $x \in G$ is 
$x + \Pb = (x + P_i)_{i\le m}$. 
We say that $\Lb$ is {\em locally finite} if for any compact set $K$ in $G$, 
$K \cap \Lb $ is finite (equivalently each $\Lam_i$ is discrete and closed).
A {\em Delone} set is a relatively dense and uniformly discrete subset in $G$
and we say that $\Lb$ is {\em Delone} if each component $\Lam_i$ is Delone. 
Throughout this paper, all multi-sets that we consider will be assumed to be 
locally finite.

\begin{defi} \label{def-Meyer}
{\em We say that a subset $\Lambda$ of $G$ is {\em Meyer} if it is a Delone set in $G$ and satisfies $\Lambda - \Lambda \subset \Lambda + J$ for some finite set $J \in G$. Similarly we say that $\Lb$ is {\em Meyer} if each component $\Lam_i$ is Meyer.} \footnote{All model sets are Meyer sets, but the reverse is far from true. Unlike model sets, Meyer sets have many characterizations \cite{RVM0}. Another definition, which is
more commonly found in the literature, is the characterization of Jeffrey Lagarias: $S$ is Meyer if it is Delone and $S-S$ is uniformly discrete \cite{Lag}. The proof in \cite{Lag} is for real spaces. A proof for compactly generated locally compact Abelian groups can be found in \cite{BLM}. Outside this setting the two concepts are different.}
\end{defi}

\begin{defi} \label{def-flc}
{\em A locally finite multi-set $\Lb$ has {\em finite local complexity
(FLC)} if for every compact set $K \subset G$ there exists a finite set $Y\subset \supp(\Lb) :=
\bigcup_{i=1}^m \Lam_i$  such that 
$$
\forall x\in \supp(\Lb),\ \exists\, y\in Y:\ 
\Lb \cap (x+K) =  (x-y) + (\Lb \cap (y+K) ) .
$$}
\end{defi}

\begin{defi} \label{def-repetitive}
{\em A locally finite multi-set $\Lb$ is {\em repetitive} if for every compact set 
$K \subset G$, $\{t \in G : \Lb \cap K = (t + \Lb) \cap K\}$ is relatively dense; i.e.
there exists a compact set $B_K \in G$ such that every translate of
$B_K$ contains at least one element of $\{t \in G : \Lb \cap K = (t + \Lb) \cap K\}$.}
\end{defi}

In order to speak about densities and autocorrelations, we need a way to specify averages.
For this purpose we fix, once and for all, an averaging sequence $\Ak = \{A_n\}_{n \in \Nat}$ satisfying \\
(i) each $A_n$ is a compact set of $G$;\\
(ii) for all $n$, $A_n \subset {A}^{\circ}_{n+1}$;\\
(iii) $\bigcup_{n \in \Nat} A_n = G$ ;\\
(iv) (the {\em van Hove} property) for all compact sets $K \subset G$, 
\[ \lim_{n \rightarrow \infty} \supreme {\frac{\theta(\partial^K(A_n))}{\theta(A_n)}} = 0,
\]
where $\partial^K(A_n) := ((K + A_n) \backslash {A^{\circ}_n}) \cup ((-K + \overline{G \backslash A_n}) \cap A_n)$ is the $K$-boundary of 
$A_n$ and $\theta$ is a Haar measure of $G$.

\medskip

\begin{defi} \label{cut-and-project1}
{\em A {\em {cut and project scheme}} (CPS) consists of a collection of spaces and mappings as follows;
\be
\begin{array}{ccccc} \label{cut-and-project-scheme1}
 G & \stackrel{\pi_{1}}{\longleftarrow} & G \times H & \stackrel{\pi_{2}}
{\longrightarrow} & H \\ 
 && \bigcup \\
 && \widetilde{L}
\end{array} 
\ee
where $G$ is a compactly generated locally compact Abelian group, $H$ is a locally 
compact Abelian group, $\pi_{1}$ and $\pi_{2}$ are the canonical projections,
$ \widetilde{L} \subset {G \times H}$ is a lattice,  i.e.\  
a discrete subgroup for which the quotient group 
$(G \times H) / \widetilde{L}$ is
compact, $\pi_{1}|_{ \widetilde{L}}$ is injective, 
and $\pi_{2}(\widetilde{L})$ is dense in $H$.}
\end{defi}

For a subset $V \subset H$, we denote $\Lam(V) := \{ \pi_{1}(x) \in G : x \in \widetilde{L}, \pi_{2}(x) \in V\}$. 

\begin{defi} \label{modelSetDef}
{\em A {\em model set} in $G$ is a subset $\Gam$ of $G$ for which, up to 
translation, $\Lam(W^{\circ}) \subset \Gam \subset \Lam(W)$, $W$ is compact in $H$, 
$W = \overline{W^{\circ}} \neq \emptyset$.
The model set $\Gam$ is {\em regular} if the boundary  
$\partial{W} = W \backslash W^{\circ}$ of $W$
is of (Haar) measure $0$.
We say that $\Gb$ is a {\em model multi-set} (resp. {\em regular model multi-set}) if each $\Gam_i$ is a model set (resp. regular model set) with respect to the same CPS. }
\end{defi}

One should note here that as $\pi_2$ need not be $1-1$ on $\widetilde{L}$, the model
set $\Gam$ need not actually be of the form $\Lam(V)$ for any set $V \subset H$. Nonetheless it is hemmed in between two such sets differing only by points on the boundary of the window $W$. The reason for this slight
relaxation in the definition of model sets is to allow for the assumption of irredundancy (\ref{irredundCond})
later in the paper.

In talking about a model set there is always an implied CPS from which it arises. When we need to be more precise we explicitly mention the CPS.

\section{Two dynamical hulls} \label{2DH}
In this section we introduce two dynamical hulls; one with a topology based on local structure and the other with a topology based on the average overall structure.

For the rest of the paper, 
\smallskip

\noindent
\fbox{$\Lb$ is a locally finite $m$-coloured multi-set for which each $\Lam_i$ is 
relatively dense.}

\smallskip
\smallskip

Define $L$
to be the subgroup of $G$ generated by all the sets $\Lam_i: i\leq m$. This group is countable. In
fact since $G$ is compactly generated it is also $\sigma$-compact: $G = \bigcup_{n=1}^\infty K_n$
where each $K_n$ is compact and $K_n\subset K_{n+1}$ for all $n$. Now $\Lam_i$ is locally finite, so $\Lam_i \cap K_n$ is finite for all
$n$. Thus $\langle \Lam_i \cap K_n \rangle$ is countable, whence also $\langle \Lam_i\rangle$
and $L$.

Let $\widetilde{\mathcal{D}} = \widetilde{\mathcal{D}(m)} $ be the set of all locally finite $m$-multi-sets in $G$. We set, for each compact set $K \subset G$ and neighbourhood $V$ of $0$ in $G$,
\be \label{uniformityFor-X}
 U_{K,V} := \{(\Lb', \Lb'') \in \widetilde{\mathcal{D}} \times \widetilde{\mathcal{D}} : 
(s+ \Lb') \cap K =  \Lb'' \cap K ~\mbox{for some}~ s \in V \}.
\ee
The set of $U_{K,V}$ forms a fundamental set of entourages of a uniformity on $\widetilde{\mathcal{D}}$ for which 
the sets of the form $U_{K,V}[\Lb''] := U_{K,V} \cap (\widetilde{\mathcal{D}} \times \{\Lb''\})$ form
a basis for the neighbourhoods of $\Lb'' \in \widetilde{\mathcal{D}}$. We call this topology the {\em local topology}. The set $\widetilde{\mathcal{D}}$ is complete with respect to the uniformity (see \cite{martin}) and is a Hausdorff space. The local hull $\mathbb{X}(\Lb)$ is defined as the closure of the $G$-orbit of $\Lb$ in $\widetilde{\mathcal{D}}$, i.e. $\mathbb{X}(\Lb) := \overline{\{ g + \Lb : g \in G\}}$.

\medskip

Now we construct the autocorrelation group $\mathbb{A}(\Lb)$. 
Let $\Lb', \Lb''$ be two locally finite $m$-multi-sets in $\widetilde{\mathcal{D}}$.
We define 
\begin{equation} \label{pseudoMetric}
d(\Lb', \Lb'') := \lim_{n \rightarrow \infty} \supreme \frac{\sum_{i=1}^m \sharp((\Lam_i' \, \triangle \, \Lam_i'') \cap A_n)}{\theta(A_n)}.
\end{equation}
Here $\triangle$ is the symmetric difference operator and $\{A_n\}$ is our averaging sequence.
This pseudo-metric is $G$-invariant.

For each open neighbourhood $V$ of $0$ in $G$ and each $\epsilon > 0$, define 
\be U(V, \epsilon) := \{ (x,y) \in G \times G : d(-v+x+\Lb, y+\Lb) < \epsilon ~\mbox{for some}~ v \in V\}. 
 \label{acUniformity}
\ee
Let $\mathcal{U} = \{U(V, \epsilon) \subset G \times G : V ~\mbox{is an open neighbourhood of}~ 0~ \mbox{in} ~G ~ \mbox{and}~ 
\epsilon > 0\}$. Then $\mathcal{U}$ forms a fundamental set of entourages for a uniformity on $G$.
Since each $U(V, \epsilon)$ is $G$-invariant, we obtain a topological group structure on $G$. We will call it the {\em autocorrelation topology}. Let $\mathbb{A}(\Lb)$ be the Hausdorff completion of $G$ in this topology, which is a new topological group (see \cite[Chapter III, \S 3.4]{Bour}). 

The meaning of $\mathbb{A}(\Lb)$ can be described as follows. The pseudo-metric $d$ is
defined on $\widetilde{\mathcal{D}}$. Define an equivalence relation on $\widetilde{\mathcal{D}}$  by
$\Lb' \equiv \Lb'' \Leftrightarrow d(\Lb', \Lb'')=0$. Let $\mathcal{D} := \widetilde{\mathcal{D}}/
\equiv$ and let $d$ also denote the resulting $G$-invariant metric on $\mathcal{D}$. Then $\mathcal{D}$ is a complete space \cite{MS}
and its elements can  be identified with actual multi-sets in $G$ up to density $0$ changes.
Now following the idea of the local topology we can put a new uniformity on $\mathcal{D}$, which mixes
the $d$-topology with the standard topology of $G$, by
using the sets
$$U'(V,\epsilon) = \{(\Lb', \Lb'') \in \mathcal{D} \times \mathcal{D} : d(-v + \Lb', \Lb'') < \epsilon \; {\textrm{for some}} \; v \in V \} \, .$$
Then $\mathcal{D}$ is complete in this topology too. Just as in forming the local hull, we can start with $\Lb$ and form the closure $\A'$ of the orbit $G + \Lb$ in $\mathcal{D}$. The mapping
$G \rightarrow \mathcal{D}$, $t \mapsto t+\Lb$  allows us to pull the topology back to $G$. This is
the content of (\ref{acUniformity}). The $\mathbb{A}(\Lb)$ is then $G$-homeomorphic with hull $\A'$
through the continuous extension of $t \mapsto t+\Lb$. For this reason we may refer to $\mathbb{A}(\Lb)$
as the autocorrelation hull.

For $y \in G$ and $U \in \mathcal{U}$, define $U[y] = \{x \in G \,: \, (x,y) \in U\}$. For each $\epsilon > 0$ let $P_{\epsilon} = \{x \in G \, : \, d(x+\Lb, \Lb) < \epsilon\}$ . Then 
$U(V, \epsilon)[0] = P_{\epsilon} + V$. Note that for any $\epsilon > 2 d(\Lb, \emptyset)$, $P_{\epsilon} = G$.

The following result is a generalization of a result found in \cite{MS}.
\begin{prop} \label{A-and-P}
Let $\Lb$ be a locally finite multi-set in $G$.
Then $\mathbb{A}(\Lb)$ is compact if and only if for all $\epsilon > 0$, $P_{\epsilon}$ is relatively dense in $G$.
\end{prop}

\noindent {\sc Proof.}  Suppose that $\mathbb{A}(\Lb)$ is compact. Since $\mathbb{A}(\Lb)$ is the completion of $G$, $G$ is precompact. So for any $\epsilon > 0$ and open neighbourhood $V$ of $0$ in $G$ whose closure is compact, there are $t_j \in G$ with $1 \leq j \leq M$ such that 
\[G \subset \bigcup_{j=1}^M (t_j +U(V,\epsilon)[0]) = \bigcup_{j=1}^M(t_j + P_{\epsilon} + V) \subset P_{\epsilon} + K,
\] where $K := \bigcup_{j=1}^M \overline{(t_j + V)}$ is compact.
Therefore $P_{\epsilon}$ is relatively dense for all $\epsilon > 0$.

Conversely, we assume that $P_{\epsilon}$ is relatively dense for all $\epsilon > 0$.
Let $\epsilon > 0$ and $V'$ be an open neighbourhood of $0$ in $G$.
From the assumption, $G \subset P_{\epsilon} + K'$ for some compact set $K'$. We can cover $K'$ with finite translations of $V'$ i.e. there are $t_1,\dots,t_L \in G$ such that $K' \subset \bigcup_{j=1}^L (t_j + V')$. Thus 
\[G \subset P_{\epsilon} + K' \subset \bigcup_{j=1}^L (t_j + P_{\epsilon} + V') = 
\bigcup_{j=1}^L U(V', \epsilon)[t_j].
\] Hence $G$ is precompact and $\mathbb{A}(\Lb)$ is compact. \qed

\medskip
 
Let us assume that $P_{\epsilon}$ is relatively dense for all $\epsilon > 0$. Recall the subgroup $L$
generated by all the points $\Lam_i, i \le m$.
It is immediate that $P_{\epsilon} \subset L$ for any $\epsilon < 2 d(\Lb, \emptyset)$.
We can define a uniformity on $L$ using the pseudo-metric $d$. Of course the $\{P_\epsilon\}$ form
a basis for a fundamental system of neighbourhoods of $0$ in the corresponding topology. We define $H$ to be a (Hausdorff) completion of $L$ in this uniformity. Then $H$ is a locally compact Abelian group
(see \cite{BM}). By definition there exists a uniformly continuous mapping $\phi : L \rightarrow H$ such that 
$\phi(L)$ is dense in $H$. 
Now we can construct a cut and project scheme:
\be \label{CPS1}
\begin{array}{ccccc}
 G & \stackrel{\pi_{1}}{\longleftarrow} & G \times H & \stackrel{\pi_{2}}
{\longrightarrow} & H \\ 
  && \cup \\
 L & \longleftarrow & \widetilde{L} & \longrightarrow & \phi(L)\\
  &&  \\
x & \longleftarrow & (x, \phi(x)) & \longrightarrow & \phi(x),
\end{array} 
\ee
where $\widetilde{L} = \{(x, \phi(x)) : x \in L\}$.
Here $\widetilde{L}$ is relatively dense and a discrete subgroup in $G \times H$, and so the factor group $\mathbb{T}(\Lb) := (G \times H)/ \widetilde{L}$ is compact (see \cite{BM}). The pseudo-metric $d$ on $L$ determines a corresponding metric $d^H$ on $H$. Let $B_{\epsilon}^H$ denote the corresponding open ball of radius $\epsilon$ in $H$. Then it is a standard feature of completions that 
$\phi(P_{\epsilon}) = \phi(L) \cap B_{\epsilon}^H$ \cite{Bour}. We define 
\[\iota : G \rightarrow (G \times \phi(L))/\widetilde{L} ~\mbox{by}~ \iota(x) = (x,0) + \widetilde{L}.
\]

\begin{prop} \label{A-and-T}
Let $\Lb$ be a locally finite multi-set in $G$. Suppose that $P_{\epsilon}$ is relatively dense for all $\epsilon > 0$. Then
$\mathbb{A}(\Lb) \cong \mathbb{T}(\Lb)$.
\end{prop}

\noindent {\sc Proof.} $\T(\Lb)$ may be viewed as the completion of $G$ under the uniform topology that is the coarsest topology on $G$ for which the map $\iota : G \longrightarrow (G \times \phi(L))/\widetilde{L}$ is continuous. For an open neighbourhood $(V \times \phi(P_{\epsilon}) + \widetilde{L})$ of $0$ in $(G \times \phi(L))/\widetilde{L}$ with $\epsilon < 2d(\Lb, \emptyset)$,
\be
\lefteqn{V \times \phi(P_{\epsilon}) + \widetilde{L} = (V - P_{\epsilon}) \times \{0\} + \widetilde{L}} \nonumber \\ 
&=& 
(V + P_{\epsilon})\times \{0\} + \widetilde{L} = \iota(V + P_{\epsilon}) = \iota(U(V, \epsilon)[0]).
\nonumber
\ee
This is the same topology on $G$ as the autocorrelation topology on $G$. Thus $\mathbb{A}(\Lb) \cong \mathbb{T}(\Lb)$.
\qed

\begin{cor}
Let $\Lb$ be a multi-set with FLC.
Suppose that there exists a continuous $G$-map $\beta : \mathbb{X}(\Lb) \rightarrow \mathbb{A}(\Lb)$. Then $P_{\epsilon}$ is relatively dense in $G$ for all $\epsilon > 0$ and
$\mathbb{A}(\Lb) \cong \mathbb{T}(\Lb)$. 
\end{cor}

\noindent {\sc Proof.} 
Since $\Lb$ has FLC, $\mathbb{X}(\Lb)$ is compact. Note that $\beta(\mathbb{X}(\Lb))$ is dense in $\mathbb{A}(\Lb)$ because 
$\beta$ is $G$-map. The continuity of $\beta$ gives us that $\mathbb{A}(\Lb)$ is compact. Thus $P_{\epsilon}$ is relatively dense in $G$ for all $\epsilon > 0$. By Proposition \ref{A-and-T} we can conclude that 
$\mathbb{A}(\Lb) \cong \mathbb{T}(\Lb)$.  \qed

\medskip
Thus we may identify $\mathbb{A}(\Lb)$ with $\mathbb{T}(\Lb)$ when it is convenient.

\begin{prop} \label{delta-compact}
Let $\Lb$ be a Meyer multi-set in $G$ and set $\Delta_i := \Lam_i - \Lam_i$ for each $i \leq m$.
If $P_{\epsilon}$ is relatively dense for all $\epsilon > 0$ then, for each $i \leq m$,
$\Delta_i$ is precompact in $L$ with the pseudo-metric topology and $\overline{\phi(\Delta_i)}$ is compact in $H$ .
\end{prop}

\noindent {\sc Proof.} Choose any $\epsilon >0$. Since each $\Lam_i$
is a Meyer set, $\Lb$ is locally finite and there exist finite sets $J_i$
with $\Lam_i -\Lam_i \subset \Lam_i + J_i$. It is easy to see that $\Delta_i - \Delta_i \subset \Delta_i + F_i$ for some finite set $F_i$ and $\Delta_i + F_i \subset \Lam_i + J_i + F_i$ for each $i \leq m$. Note that $\Lam_i + J_i + F_i$ is locally finite. From the assumption, for any given $P_{\epsilon}$ we can find a compact set $K_{\epsilon}$ such that $G \subset P_{\epsilon} + K_{\epsilon}$. So $\Delta_i \subset P_{\epsilon} + K_{\epsilon}$ for each $i \leq m$. Since we are assuming that each $\Lam_i$ is relatively dense, for small enough $\epsilon > 0$, $P_{\epsilon} \subset \cap_{i=1}^m \Delta_i$.  Then by the assumption 
\[ (\Delta_i - P_{\epsilon}) \cap K_{\epsilon} \subset (\Delta_i - \Delta_i) \cap K_{\epsilon} 
\subset (\Delta_i + F_i) \cap K_{\epsilon} = N_i, 
\]
where $N_i$ is a finite subset of $L$.
Thus $\Delta_i \subset P_{\epsilon} + N_i$ for each $i \leq m$.
So $\Delta_i$ is precompact (totally bounded) in $L$ with the pseudo-metric topology and $\overline{\phi(\Delta_i)}$ is compact  in $H$ for each $i \leq m$. \qed

\medskip

\section{Torus parametrizations for model multisets} \label{Map-on-dyn-hulls}

A torus parametrization $\mathbb{X}(\Lb)$ is a continuous $G$-map $\beta : \mathbb{X}(\Lb) \rightarrow \mathbb{A}(\Lb)$. \footnote{The terminology arises from the model set cases first studied in which
(in the set-up that we have here) $\mathbb{A}(\Lb)$ would have been a torus.}  An element $\Gb \in \mathbb{X}(\Lb)$ is {\em non-singular} for this parametrization if $\beta^{-1}(\beta(\{\Gb\}))= \{\Gb\}$. The set of non-singular elements of $\mathbb{X}(\Lb)$ consists of $G$-orbits.
In this section we show that the existence of a torus parmetrization together with the existence of a non-singular element in $\mathbb{X}(\Lb)$ ensures that $\Lb$ is a model multi-set.

\begin{prop} \label{one-ele-on-fiber}
Let $\Lb$ be a multi-set in $G$ with finite local complexity (FLC).
Suppose that  there exists a continuous $G$-map $\beta : \mathbb{X}(\Lb) \rightarrow \mathbb{A}(\Lb)$. 
If $\Gb$ is non-singular in $\mathbb{X}(\Lb)$, then given any $M \in \Z_+$, there is 
$\epsilon = \epsilon(M) > 0$ such that for all $t \in P_{\epsilon}$, $(t + \Gb) \cap A_M = \Gb \cap A_M$.
\end{prop}

\noindent {\sc Proof.} Suppose that there is a positive integer $M$ such that for every $n > 0$ there exists 
$t_n \in P_{(1/{2^n})}$ for which 
\[(t_n + \Gb) \cap A_M \neq \Gb \cap A_M.
\] 
Since $\mathbb{X}(\Gb) \subset \mathbb{X}(\Lb)$ is compact from the FLC of $\Lb$, 
$\{t_n + \Gb \}_n$ has a convergent subsequence $\{t_{n_k} + \Gb \}_k$ such that 
$t_{n_k} + \Gb \rightarrow \Gb'$. On the other hand, identifying $\mathbb{A}(\Lb)$ with 
$\mathbb{T}(\Lb)$, 
\[\beta(t_{n_k} + \Gb) = (t_{n_k}, 0) + \beta(\Gb) =(0, -\phi(t_{n_k})) + \beta(\Gb).\] 
Since $\phi(t_{n_k}) \stackrel{k \rightarrow \infty}{\longrightarrow} 0$,
$\lim_{k \rightarrow \infty} \beta(t_{n_k} + \Gb) = \beta(\Gb)$. Yet by the continuity of $\beta$, 
\[\lim_{k \rightarrow \infty} \beta(t_{n_k} + \Gb) = \beta(\Gb').\]
Then $\beta(\Gb) = \beta(\Gb')$, and 
by the assumption that $\Gb$ is non-singular in $\mathbb{X}(\Lb)$ $\Gb = \Gb'$. It follows that $t_{n_k} + \Gb \stackrel{k \rightarrow \infty}{\longrightarrow} \Gb$. This contradicts the choice of the original sequence $\{t_{n}\}$. \qed

\begin{prop} \label{sing-open-win}
Let $\Lb$ be a multi-set in $G$ with FLC.
Suppose that there exists a continuous $G$-map $\beta : \mathbb{X}(\Lb) \rightarrow \mathbb{A}(\Lb)$.
If $\Gb$ is non-singular in $\mathbb{X}(\Lb)$ then there exist $s \in G$ and non-empty open sets $ U_i \subset H$ so that  
\[\Gam_i = -s + \Lam(U_i) ~~\mbox{for each}~i \le m\]
 with respect to CPS (\ref{CPS1}).
Furthermore each $\overline{U_i}$ is compact if and only if $\Gb$ is Meyer.
\end{prop} 

\noindent {\sc Proof.} 
Since $\Gb \in \mathbb{X}(\Lb)$, we can choose $s \in G$ and a compact set $K$ for which 
$(s + \Gb) \cap K = \Lb \cap K$ and each $(s + \Gam_i) \cap K \neq \emptyset$. 
Then $s + \Gam_i \subset \Lam_i + L \subset L$ for all $i \le m$. Furthermore, since $\beta$ is 
$G$-map and $\Gb$ is non-singular, $s + \Gb$ is also non-singular. So we can translate and assume at the outset that 
$\supp(\Gb) \subset L$. 

Let $i \leq m$ and $x \in \Gam_i$. For any $M > 0$ with $x \in \Gam_i \cap A_M$, there is 
$\epsilon_x = \epsilon(M) > 0$ so that for any $y \in P_{\epsilon_x}$, $(y + \Gam_i) \cap A_M = \Gam_i \cap A_M$ by Proposition \ref{one-ele-on-fiber}. This implies that $x - y \in \Gam_i$ for any $y \in P_{\epsilon_x}$, and thus 
$x - P_{\epsilon_x} \subset \Gam_i$. Since $P_{\epsilon_x} = - P_{\epsilon_x}$, $x + P_{\epsilon_x} \subset \Gam_i$. Therefore 
\[\Gam_i = \bigcup_{x \in \Gam_i}(x + P_{\epsilon_x})~\mbox{for}~ i \leq m.\]

Recall that $\phi(P_{\epsilon_x}) = \phi(L) \cap B_{\epsilon_x}^H$ where $B_{\epsilon_x}^H$ is an open ball of radius ${\epsilon}_x$ in $H$.
So 
\[\Gam_i = \Lam (U_i) ~~\mbox{where}~ 
U_i = \bigcup_{x \in \Gam_i}(\phi(x) + B_{\epsilon_x}^H)~\mbox{is open in}~ H.\] 

From Proposition \ref{A-and-P}, $P_{\epsilon}$ is relatively dense for all $\epsilon > 0$. If $\Gb$ is Meyer, $\overline{U_i}$ is compact from Proposition \ref{delta-compact}. On the other hand, if each $\overline{U_i}$ is compact then $\Gam_i - \Gam_i \subset \Lam(\overline{U_i} - \overline{U_i})$ which is a Delone set and so $\Gam_i$ is a Meyer set.
\qed

\begin{prop} \label{non-sin-no-boundary}
Let $\Lb$ be a multi-set in $G$ with FLC.
Suppose that there exists a continuous $G$-map $\beta : \mathbb{X}(\Lb) \rightarrow \mathbb{A}(\Lb)$.
If $\Gb \in \mathbb{X}(\Lb)$ is non-singular and $\Gam_i = \Lam(U_i)$ with an open set $U_i \subset H$ for each $i \le m$ with respect to CPS (\ref{CPS1}),
then $\phi(L) \cap \partial{U_i} = \emptyset$ for all $i \le m$.
\end{prop}

\noindent {\sc Proof.} 
Suppose that $\phi(x) \in \phi(L) \cap \partial{U_i}$ for some $x \in L$. So 
$x \notin \Gam_i = \Lam(U_i)$.
There is a sequence $\{y_n\}_n$ in $\Gam_i$ such that $\phi(y_n) \in \phi(L) \cap U_i$ and $\lim_{n \rightarrow \infty} \phi(y_n) = \phi(x)$.
Since $\mathbb{X}(\Gb) \subset \mathbb{X}(\Lb)$ is compact, we can assume that $\{x - y_n + \Gb\}_n$ is a converging sequence in 
$\mathbb{X}(\Gb)$. Let $\Gb' := \lim_{n \rightarrow \infty}(x - y_n + \Gb)$. Note that $\Gam'_i$ contains 
$x$ and $\beta(\Gb') = \beta(\Gb)$. Since $\Gb$ is non-singular in $\mathbb{X}(\Lb)$, $\Gb' = \Gb$. 
However $x \notin \Gam_i$. This contradiction shows that $ \phi(L) \cap \partial{U_i} = \emptyset$.  \qed

\medskip

\noindent {\sc \bf Remark:} By this last Proposition it is possible to make the sets
$U_i$ a little nicer, replacing them by ${W_i}^{\circ}$ where $W_i := \overline{U_i}$,
since these differ only from the $U_i$ by boundary points, and these additional points
do not affect the model multi-set $\Gb$. This $W_i$ has the nice property that
$W_i$ is the closure of its interior, a property often assumed in defining model sets.

\medskip

The following proposition holds in general CPS.

\begin{prop} \label{between-interior-closure}
Let $\Lb$ be a multi-set in $G$ for which $\Lam({V_i}^{\circ}) \subset \Lam_i \subset \Lam(\overline{V_i})$ where $\overline{V_i}$ is compact and ${V_i}^{\circ} \neq \emptyset$ for $i \le m$, with respect to some CPS (see (\ref{cut-and-project-scheme1})). Then $\Lb$ has FLC, and for any $\Gb \in \mathbb{X}(\Lb)$ there exists $(-s, -h) \in G \times H$ so that 
\[-s + \Lam(h + {V_i}^{\circ}) \subset \Gam_i \subset -s + \Lam(h + \overline{V_i}) ~~ 
\mbox{for each}~ i \leq m. 
\]
If \, $\bigcup_{i=1}^m \Gamma_i \subset L$, then we can take $s=0$.
Furthermore, if there exists a continuous $G$-map $\beta : \mathbb{X}(\Lb) \rightarrow \mathbb{A}(\Lb)$ such that 
$\beta(\Lb) = (0,0) + \widetilde{L}$, then $\beta(\Gb) = (-s, -h) + \widetilde{L}$.
\end{prop}

\noindent {\sc Proof.} For each $\Gb \in \mathbb{X}(\Lb)$ we can choose $s \in G$ such that 
$s + \bigcup_{i=1}^m \Gam_i \subset L$. We first claim that $\Lb$ has FLC. Since  
$\Lam({V_i}^{\circ}) \subset \Lam_i \subset \Lam(\overline{V_i})$ for each $i \le m$, 
$$\Lam(\bigcup_{i \le m} {V_i}^{\circ}) \subset \bigcup_{i \le m}\Lam_i \subset \Lam(\bigcup_{i \le m} \overline{V_i}).$$ Note that 
$W := \bigcup_{i \le m} \overline{V_i}$ is compact. So for any compact set $K \subset G$, 
$(\Lam(W) - \Lam(W)) \cap K$ has a finite number of elements. This implies that 
for any $x \in \supp(\Lb)$ there are, 
up to translation, only finitely many different subsets of the form $(K + x) \cap \Lb$.
Therefore the claim follows.
Then $\X(\Lb)$ is compact and we can find $\{t_n\}_n \subset L$ such that $\{t_n + \Lb\}$ 
converges to $s + \Gb$ in $\mathbb{X}(\Lb)$. There exists $n_0 \in \Z_+$ such that for any $k,l \geq n_0$, $(t_k + \Lb) \cap (t_l + \Lb) \neq \emptyset$ and so $t_k - t_l \in \Lam_i - \Lam_i$ for some $i \leq m$. Thus $\phi(t_k) - \phi(t_l) \in \overline{V_i} - \overline{V_i}$. Since $\overline{V_i} - \overline{V_i}$ is compact, we can find a convergent subsequence of $\{\phi(t_n)\}_n$. Without loss of generality, we can assume that 
\[\lim_{n \rightarrow \infty} \phi(t_n) =: h \in H.\]
For each $i \le m$, if $z \in \Lam(h + {V_i}^{\circ})$, then there is $n_{i,1} \in \Z_+$ such that $\phi(z) - \phi(t_n) \in {V_i}^{\circ}$ for any $n \ge n_{i,1}$. So $z - t_n \in 
\Lam({V_i}^{\circ}) \subset \Lam_i$ i.e. 
$z \in t_n + \Lam_i$. Thus $z \in s + \Gam_i$. 
This implies 
\[\Lam(h + {V_i}^{\circ}) \subset s + \Gam_i ~~ \mbox{for each}~ i \leq m.\] 
On the other hand, if $z \in s + \Gam_i$, then $z \in t_n + \Lam_i$ for large $n$.
So $\phi(z) \in \phi(t_n) + \overline{V_i}$ for large $n$ and $\phi(z) \in h + \overline{V_i}$. Thus $z \in \Lam(h + \overline{V_i})$. Therefore 
\[\Lam(h + {V_i}^{\circ}) \subset s + \Gam_i \subset \Lam(h + \overline{V_i})~~ \mbox{for each}~ 
i \leq m.\] 
Note that $\beta(t_n + \Lb) = \iota(t_n) + \beta(\Lb)$ and 
$\iota(t_n) = (t_n, 0) + \widetilde{L} = (0, -\phi(t_n)) + \widetilde{L}$. 
Thus 
\begin{eqnarray*}
\beta(\Gb) &=& \iota(-s) + \beta(s + \Gb) = \iota(-s) + \lim_{n \to \infty} \beta(t_n + \Lb) \\
& =& \iota(-s) + \lim_{n \to \infty} \iota(t_n) + \beta(\Lb) = (-s,0) + (0,-h) + \widetilde{L} \\
&=& (-s,-h) + \widetilde{L}.
\end{eqnarray*}
 \qed

\begin{cor} \label{most-interior-win}
Let $\Lb$ be a multi-set with FLC.
Suppose that there exists a continuous $G$-map $\beta : \mathbb{X}(\Lb) \rightarrow \mathbb{A}(\Lb)$, and 
$\Lam_i = \Lam(W_i)$ with respect to CPS (\ref{CPS1}) where $W_i \subset H$ is compact and $W_i = \overline{{W_i}^{\circ}} \neq \emptyset$ for each $i \leq m$.
If $\Gb$ is non-singular in $\mathbb{X}(\Lb)$, then there is $(-s, -h) \in G \times H$ so that 
\[\Gam_i = -s + \Lam(h + {W_i}^{\circ})~~\mbox{for each}~ i \le m.
\] 
\end{cor}

\noindent {\sc Proof.} If $\Gb$ is non-singular in $\mathbb{X}(\Lb)$, 
from Proposition \ref{sing-open-win} and Proposition \ref{non-sin-no-boundary} there exists $s \in G$ such that $\Gam_i = -s + \Lam(U_i)$ where 
$U_i \subset H$ is open and $\phi(L) \cap \partial{U_i} = \emptyset$ for each $i \le m$. Then from Proposition \ref{between-interior-closure} there exists $(-s, -h) \in G \times H$ (same $s$) such that 
$-s + \Lam(h + {W_i}^{\circ}) \subset \Gam_i \subset -s + \Lam(h + W_i)$ for each 
$i \le m$. Thus for each $i \leq m$,
\[\Lam(h + {W_i}^{\circ}) \subset \Lam(U_i) \subset \Lam(h + W_i).\]
Since 
$U_i \backslash (h + W_i)$ is open and $\phi(L)$ is dense in $H$, $U_i \subset h + W_i$. 
Replacing $U_i$ by $(h + {W_i}^{\circ}) \cup  U_i$, we can assume that $h + {W_i}^{\circ} \subset U_i \subset h + W_i$.
Since $U_i$ is an open set, $U_i \subset h + {W_i}^{\circ}$. 
Therefore $U_i = h + {W_i}^{\circ}$ and 
$\Gam_i = -s + \Lam(h + {W_i}^{\circ})$. \qed

\begin{prop} \label{conti-nonemptyModel-Model}
Let $\Lb$ be a multi-set in $G$ with FLC and repetitivity. 
Suppose that there exists a continuous $G$-map $\beta : \mathbb{X}(\Lb) \rightarrow \mathbb{A}(\Lb)$ and that $\Lam({V_i}^{\circ}) \subset \Lam_i \subset \Lam(\overline{V_i})$ where $\overline{V_i}$ is compact, ${V_i}^{\circ} \neq \emptyset$, and $\partial{V_i}$ has empty interior for each $i \le m$, with respect to CPS (\ref{CPS1}).  
Then there exists a non-singular element $\Lb'$ in $\mathbb{X}(\Lb)$ such that 
$\Lam'_i = \Lam(W_i)$ where $W_i$ is compact and $W_i = \overline{{W_i}^{\circ}}$ for each $i \leq m$ with respect to the same CPS (\ref{CPS1}), and so that 
for each $\Gb \in \mathbb{X}(\Lb)$ there exists $(-s, -h) \in G \times H$ so that 
\[-s + \Lam(h + {W_i}^{\circ}) \subset \Gam_i \subset -s + \Lam(h + W_i) ~~\mbox{for each}~ i \le m.
\] 
\end{prop}

\noindent {\sc Proof.} 
From Proposition \ref{between-interior-closure}, for any $\Gb \in \mathbb{X}(\Lb)$ there exists $(-s, -h) \in G \times H$ so that 
\[-s + \Lam(h + {V_i}^{\circ}) \subset \Gam_i \subset -s + \Lam(h + \overline{V_i}) ~~ 
\mbox{for each}~ i \leq m. 
\]
By the Baire category theorem, there exists $h' \in H$ such that 
$$\phi(L) \cap (h' + \bigcup_{i \le m}\partial{V_i}) = \emptyset.$$
So we can find $\Lb' \in \mathbb{X}(\Lb)$ such that $\beta(\Lb') = (0, -h') + \widetilde{L}$ and 
$\Lam'_i = \Lam(h' + {V_i}^{\circ})$ for $i \le m$.
For each $i \le m$, let $W_i = h' + \overline{{V_i}^{\circ}}$. Then $W_i = \overline{{W_i}^{\circ}}$, $W_i$ is compact, and $\Lam'_i = \Lam(W_i)$.
Since repetitivity of $\Lb$ is equivalent to the minimality of $\mathbb{X}(\Lb)$ \cite{Fur, martin}, $\mathbb{X}(\Lb) = \mathbb{X}(\Lb')$. Also $\Lb'$ is non-singular by construction.
We conclude using  Proposition \ref{between-interior-closure} again. \qed

\begin{prop} \label{boundary-measure}
Let $\Lb$ be a multi-set with FLC.
Suppose that there exists a continuous $G$-map $\beta : \mathbb{X}(\Lb) \rightarrow \mathbb{A}(\Lb)$ which 
is one-to-one a.e. $\mathbb{A}(\Lb)$, and $\Lam_i = \Lam(W_i)$ with respect to CPS (\ref{CPS1}) where $W_i$ is compact and $W_i = \overline{{W_i}^{\circ}} \neq \emptyset$ for each $i \leq m$.
Then 
\[\theta_H(\partial{W_i}) = 0 ~~\mbox{for each}~ i \le m,\]
where  $\theta_H$ is any Haar measure in $H$.
\end{prop}

\noindent {\sc Proof.} 
From Corollary \ref{most-interior-win}, for any non-singular element $\Gb \in \mathbb{X}(\Lb)$, 
$\Gam_i = -s + \Lam(h + {W_i}^{\circ})$. By Proposition \ref{non-sin-no-boundary}, we get  
$\phi(L) \cap (h + \partial{W_i}) = \emptyset$.
So by the assumption on $\beta$, $-s + \Lam(h + \partial{W_i}) = \emptyset$ a.e. $(-s, -h) + \widetilde{L} \in \mathbb{A}(\Lb)$.
By \cite{RVM1}, there is a normalization of $\theta_H$ for which $\dens(-s + \Lam(h + \partial{W_i})) = \theta_H(\partial{W_i})$ a.e. $\mathbb{A}(\Lb)$. 
Thus $\theta_H(\partial{W_i}) = 0$ for each $i \leq m$.
\qed

\medskip
We note that if $\Lb$ is a Meyer multi-set then $\Lb$ has FLC.

\begin{theorem} 
Let $\Lb$ be a Meyer multi-set with repetitivity.
Suppose that there exists a continuous $G$-map $\beta : \mathbb{X}(\Lb) \rightarrow \mathbb{A}(\Lb)$ which 
is one-to-one a.e. $\mathbb{A}(\Lb)$.
Then there exists $\Lb' \in \mathbb{X}(\Lb)$ such that $\Lam'_i = \Lam(W_i)$ with respect to CPS (\ref{CPS1}) where $W_i = \overline{{W_i}^{\circ}} \neq \emptyset$ and 
$W_i$ is compact for each $i \leq m$, and so 
for each $\Gb \in \mathbb{X}(\Lb)$ there exists $(-s, -h) \in G \times H$ so that 
\[-s + \Lam(h + {W_i}^{\circ}) \subset \Gam_i \subset -s + \Lam(h + W_i) ~~\mbox{for each}~ i \le m.
\] 
Furthermore 
\[\theta_H(\partial{W_i}) = \emptyset ~~\mbox{for each}~ i \leq m.
\]
In other words, for each $\Gb \in \mathbb{X}(\Lb)$, $\Gb$ is a regular model multi-set.
\end{theorem}

\noindent {\sc Proof.} 
Since $\beta$ is one-to-one a.e. $\mathbb{A}(\Lb)$, there exists a non-singular element $\Lb' \in \mathbb{X}(\Lb)$. So 
from Proposition \ref{sing-open-win} and \ref{non-sin-no-boundary} we can suppose that $\Lam'_i = \Lam(W_i)$, 
$W_i = \overline{{W_i}^{\circ}} \neq \emptyset$, and $W_i$ is compact for each $i \le m$. Note that $\mathbb{X}(\Lb') = \mathbb{X}(\Lb)$ by the repetitivity of $\Lb$ (see \cite{Fur, martin}).
Thus applying Proposition \ref{between-interior-closure} and Proposition \ref{boundary-measure} we can conclude that for each $\Gb \in \mathbb{X}(\Lb)$, $\Gb$ is a regular model multi-set. \qed

\section{Torus parameterizations from model multi-sets}

\noindent
We consider a cut and project scheme :
\be \label{cut-and-project2}
\begin{array}{ccccc}
 G & \stackrel{\pi_{1}}{\longleftarrow} & G \times H & \stackrel{\pi_{2}}
{\longrightarrow} & H \\ 
  && \cup \\
  & & \widetilde{L} & & 
\end{array} 
\ee
where $H$ is a locally compact Abelian group, $\pi_1$ and $\pi_2$ are canonical maps, $\widetilde{L}$ is a lattice in $G \times H$, 
$\pi_1|_{\widetilde{L}}$ is one-to-one, and $\pi_2(\widetilde{L})$ is dense in $H$.
Let $L = \pi_1(\widetilde{L})$. We define $\phi : L \rightarrow H$ by 
$\phi(x) = \pi_2(\pi_1^{-1}(x))$.

Suppose that $\Lb$ is a multi-set in $G$ and that for each $i\le m$,
$\Lam({W_i}^{\circ}) \subset \Lam_i \subset \Lam(W_i)$, $W_i$ is compact in $H$, and $W_i = \overline{{W_i}^{\circ}} \neq \emptyset$ with respect to CPS (\ref{cut-and-project2}). 

There are two things that we can do to tighten up the cut and project scheme without altering
the multiset $\Lb$. 

First we can require that the group $H_0 :=\langle W_i : i\le m\rangle$
is all of $H$. To see this, note first that the subgroup of $H$ generated by the ${W_i}^{\circ}$
is an open, hence closed subgroup of $H$, and so it contains all the $W_i$ and must then 
be $H_0$. So $H_0$ is closed. Let $\widetilde{L_0} :=\widetilde{L} \cap (G\times H_0)$ and
$L_0$ its $\pi_1$ projection. Now we can check that we can replace $H, \widetilde{L}$
by $H_0, \widetilde{L_0}$ in the cut and project scheme to get another one with the desired property.

Once we know that $H$ is generated by the windows, then it follows that $\langle \Lam_i : i\le m \rangle = L$. In fact $\phi(\Lam_i)$ is dense in $W_i$, so $\phi(\langle \Lam_i : i\le m \rangle)$
is dense in $H$. Then for any coset $x + \langle \Lam_i : i\le m \rangle \subset L$,
$\phi(x) +  \phi(\langle \Lam_i : i\le m \rangle)$ is also dense in $H$ and so has a point
$\phi(x) + \phi(u) \in W_1^\circ$, say, where $u \in \langle \Lam_i : i\le m \rangle$. Since $x+u \in L$,
it is in fact in $\Lam_1$ and so 
we obtain $x \in \langle \Lam_i : i\le m \rangle$.

A consequence of this is that, as we saw at the start of Sec.\,\ref{2DH}, $L$, and so also $\widetilde{L}$, is now countable.

Let $I := \{t \in H : t + W_i = W_i ~\mbox{for all}~ i \leq m\}$. Translations $t$ in $H$ of this form indicate certain redundancy in $H$. The second way in which we can tighten up the CPS
is to remove this redundancy in $H$ by factoring out the subgroup $I$.
We define $H' := H/I$, $\psi : L \rightarrow H'$ by $\psi(x) = \phi(x) + I$, and $\widetilde{L'} := \{(x, \psi(x)) \in G \times H' : x \in L \}$. Then $\widetilde{L'}$ is a lattice in $G \times H'$, i.e. $\widetilde{L'}$ is a discrete subgroup for which $(G \times H')/\widetilde{L'}$ is compact. Note that $W_i + I = W_i$ and ${W_i}^{\circ} + I = {W_i}^{\circ}$ for all $i \leq m$. Thus for all $i \leq m$
\[\Lam(c + W_i + I) = \Lam(c + W_i) ~~\mbox{and}~
\Lam(c + {W_i}^{\circ} + I) = \Lam(c + {W_i}^{\circ})  ~~\mbox{for any}~ c \in H.\] 
Let $W'_i$ denote $W_i + I$ in $H'$.
Then we can construct a new cut and project scheme :
\be \label{cut-and-project3}
\begin{array}{ccccc}
 G & \stackrel{\pi'_{1}}{\longleftarrow} & G \times H' & \stackrel{\pi'_{2}} {\longrightarrow} & H' \\ 
  && \cup &&\\
 L & \longleftarrow & \widetilde{L'} & \longrightarrow & \psi(L) \\
  & & & & \\
x & \longleftarrow & (x, \psi(x)) & \longrightarrow & \psi(x) \,
\end{array} 
\ee
and we get $\Lam({W'_i}^{\circ}) \subset \Lam_i \subset \Lam(W'_i)$, $W'_i$ is compact in $H'$, 
$W'_i = \overline{{W'_i}^{\circ}} \neq \emptyset$ for $i \leq m$, $\{t \in H' : t + W'_i = W'_i ~\mbox{for all}~ i \leq m\} = \{0\}$, with respect to CPS (\ref{cut-and-project3}). 
Furthermore since $I$ is a closed subgroup of $H$, if $\theta_{H}({\partial{W_i}}) = 0$ where $\theta_{H}$ is a Haar measure in $H$, then 
$\theta_{H'}({\partial{W'_i}}) = 0$ where $\theta_{H'}$ is a Haar measure in $H'$ (see \cite[Theorem 3.3.28]{RS}).

Thus without loss of generality we will assume both that $H$ is generated by the windows of $\Lb$ and
\begin{equation} \label{irredundCond}
 \{t \in H : t + W_i = W_i ~\mbox{for all}~ i \leq m\} = \{0\} \, , 
 \end{equation}
a situation that we will refer to as {\em irredundancy}. All subsequent CPS will be assumed reduced into this form.

\medskip

In this section we establish the existence of the `torus parametrization'. The following theorem is proved after a sequence of auxiliary propositions.

\begin{theorem} \label{main}
Let $\Lb$ be a multi-set in $G$. Suppose that $\Lam({W_i}^{\circ}) \subset \Lam_i \subset \Lam(W_i)$, $W_i$ is compact, $W_i = \overline{{W_i}^{\circ}} \neq \emptyset$, and 
$\theta_{H}({\partial{W_i}}) = 0$ for all $i \leq m$ with respect to some irredundant CPS. Then there is a continuous $G$-map $\beta : \mathbb{X}(\Lb) \rightarrow \mathbb{A}(\Lb)$ which is one-to-one a.e. $\mathbb{A}(\Lb)$.
\end{theorem}

\begin{prop} \label{uniquely-chosen-c}
Let $\Lb$ be a multi-set in $G$. Suppose that $\Lam({W_i}^{\circ}) \subset \Lam_i \subset \Lam(W_i)$,
$W_i$ is compact, and $W_i = \overline{{W_i}^{\circ}} \neq \emptyset$ for all $i \le m$ with respect to 
some irredundant CPS. 
Then for any $\Gb \in \mathbb{X}(\Lb)$ with $\bigcup_{i = 1}^m \Gam_i \subset L$, 
$$\bigcap \{\phi(t) - W_i : t \in \Gam_i, i \leq m\}$$ contains exactly one element $c_{\Gbt}$ in $H$ and 
$\overline{\phi(\Gam_i)} = c_{\Gbt} + W_i$ for each $i \leq m$. Furthermore, $c_{\Lbt} = 0$.
\end{prop}

\noindent {\sc Proof.} Let $\Gb \in \X(\Lb)$ with $\cup_{i = 1}^m \Gam_i \subset L$. We claim that \[\bigcap \{\phi(t)- W_i : t \in \Gam_i, i \leq m\} \neq \emptyset.\] 
Suppose that $ \bigcap \{\phi(t) - W_i : t \in \Gam_i, i \leq m\} = \emptyset $. Since each $\phi(t) - W_i$ is compact, there exists a compact $K \subset G$ such that $ \bigcap \{\phi(t) - W_i : t \in \Gam_i \cap K, i \leq m\} = \emptyset$ and $\Gb \cap K \neq \emptyset$. 
Since $\Gb \in \mathbb{X}(\Lb)$ and $\bigcup_{i=1}^m \Gam_i \subset L$,
we can find $t_0 \in L$ such that $\Gb \cap K = (t_0 + \Lb) \cap K$. 
Then for any $i \leq m$ and any  $t \in \Gam_i \cap K$, $t \in t_0 + \Lam_i$ and $\phi(t_0) \in \phi(t) - W_i$. This is a contradiction.

For any $c \in \bigcap \{\phi(t)- W_i : t \in \Gam_i, i \leq m\}$, $\phi(\Gam_i) \subset c + W_i$ for all $i \leq m$.
Since $W_i = \overline{{(W_i)}^{\circ}}$, there is $h \in H$ such that $\overline{\phi(\Gam_i)} = h + W_i$ for all $i \leq m$ by Proposition \ref{between-interior-closure}. 
So $h + W_i \subset c + W_i$.
We claim that $h + W_i = c + W_i$. In fact, $W_i \subset (c - h) + W_i$. Let $x = c - h$.
Since $W_i - W_i$ is compact, $\overline{\{nx : n \in \Z_+\}}$ is compact. For any neighbourhood 
$V$ of $0$ in $H$, $\{V + nx : n \in \Z_+\}$ is an open cover of $\overline{\{nx : n \in \Z_+\}}$,
so there are $n_1, \dots, n_k \in \Z_+$ such that $\{V + n_j x : 1 \le j \le k\}$ covers 
$\{nx : n \in \Z_+\}$. For $n \in \Z_+$ with $n > \max\{n_1, \dots, n_k\}$, $nx \in  n_j x + V$ for some $1 \le j \le k$ and so $V$ contains $p x$ where  $p := n - n_j \in \Z_+$.
So $W_i \subset x + W_i \subset 2x + W_i \subset \dots \subset p x + W_i \subset V+ W_i $. Since $V$ is arbitrary, $W_i = x + W_i = c - h + W_i$ for $i \le m$. Thus $c = h$, since $0 = \{t \in H : t + W_i = W_i, i \leq m\}$. So there is a unique $c_{\Gbt} \in H$ such that 
\begin{equation} \label{def of c}
\{c_{\Gbt}\} = \bigcap \{\phi(t)-W_i : t \in \Gam_i, i \leq m\}.
\end{equation} 
Since  $0 \in \bigcap \{\phi(t)-W_i : t \in \Lambda_i, i \leq m\}$, \, $c_{\Lbt} = 0$.
\qed

\medskip

As usual we will define $\mathbb{T}(\Lb) := (G \times H)/\widetilde{L}$ in the irredundant CPS.

\begin{cor} \label{def-a-map}
Let $\Lb$ be a multi-set in $G$. Suppose that $\Lam({W_i}^{\circ}) \subset \Lam_i \subset \Lam(W_i)$,
$W_i$ is compact, and $W_i = \overline{{W_i}^{\circ}} \neq \emptyset$ for all $i \le m$, with respect to some irredundant CPS. 
Then the map $\gamma : \mathbb{X}(\Lb) \rightarrow \mathbb{T}(\Lb)$ given by $\Gb \mapsto \gamma(\Gb) = (-s, -c_{s+\Gbt}) + 
\widetilde{L}$, where $s$ is any element of $G$ for which $s + \bigcup_{i=1}^m \Gam_i  \subset L$ and $c_{s+\Gbt}$ is given by (\ref{def of c}), is a well-defined G-map. Furthermore, 
for any $\Gb \in \X(\Lb)$,
\[\overline {\phi(\Gamma_i)} = -\phi(s) + c_{s + \Gbt} + W_i ~~\mbox{for all}~i \le m \] 
and $\gamma(\Lb) = (0,0) + \widetilde{L}$.
\end{cor}

\noindent {\sc Proof.} Suppose that $s \neq s_1$ where $  s +\bigcup_{i=1}^m \Gam_i \subset L$ and 
$s_1 + \bigcup_{i=1}^m \Gam_i  \subset L$. Note that $s - s_1 \in L$, and 
\be
(-s_1, -c_{s_1+\Gbt}) + \widetilde{L} & = & (-(s +l), -c_{s +l+\Gbt}) + \widetilde{L} ~~\mbox{for some}~ l \in L \nonumber \\ 
 & = & (-s -l, -\phi(l)-c_{s +\Gbt}) + \widetilde{L}  \nonumber \\
 & = & (-s, -c_{s +\Gbt}) + \widetilde{L}.
\ee
For any $g \in G$, 
\[ \gamma(g + \Gb) = (-s +g, -c_{s +\Gbt}) + \widetilde{L} = \iota(g) + \gamma(\Gb).
\]
The rest follows directly from Proposition \ref{uniquely-chosen-c}.
\qed

\begin{prop} \label{continuity}
Let $\Lb$ be a multi-set in $G$. Suppose that $\Lam({W_i}^{\circ}) \subset \Lam_i \subset \Lam(W_i)$,
$W_i$ is compact, and $W_i = \overline{{W_i}^{\circ}} \neq \emptyset$ for all $i \le m$, with respect to some irredundant CPS. 
The mapping $\gamma : \mathbb{X}(\Lb) \rightarrow \mathbb{T}(\Lb)$ defined in Corollary \ref{def-a-map} is continuous and surjective.
\end{prop}

\noindent {\sc Proof.} Let $\Gb \in \mathbb{X}(\Lb)$ and $\gamma(\Gb) = (-s, -c_{s+ \Gbt}) + \widetilde{L}$ for some $s \in G, c_{s+ \Gbt} \in H$. Let $U$ be an open neighbourhood of $0$ in 
 $G$ and let $U'$ be an open neighbourhood of $-c_{s+ \Gbt}$ in $H$. Let $U_0 = U \cap (-U)$. Since 
$c_{s+ \Gbt} = \bigcap \{\phi(t) - W_i : t \in s + \Gam_i, i \leq m\}$, 
\[\bigcap \{(-\phi(t) + W_i) \backslash U' : t \in s + \Gam_i, i \leq m\} = \emptyset.\]
Each of $(-\phi(t) + W_i) \backslash U'$ is closed and so there exists a compact $K$ in $G$ 
with $(s+\Gb)\cap K \ne \emptyset$ such that 
\[\bigcap \{(-\phi(t) + W_i) \backslash U' : t \in (s + \Gam_i) \cap K, i \leq m\} = \emptyset\] i.e. 
$\bigcap \{-\phi(t) + W_i: t \in (s + \Gam_i) \cap K, i \leq m\} \subset U'$.
For any $\Gb' \in U_{-s+K, U_0}[\Gb]$ (see (\ref{uniformityFor-X})), 
\[(r  + \Gb' ) \cap (-s + K) = \Gb \cap (-s + K) ~\mbox{for some}~ r \in U_0.
\]
So $( r +s +\Gb' ) \cap K = (S +\Gb ) \cap K$.
Then 
\[\bigcap \{-\phi(t) + W_i : t \in (r + s+\Gam'_i ) \cap K, i \le m\}  \subset U'.
\]
This shows 
\[\gamma(\Gb') \in (-r-s, U') + \widetilde{L} \subset (-s+ U_0, U') + \widetilde{L} \subset 
(-s + U, U') + \widetilde{L}.
\]
Therefore $\gamma$ is continuous.
Furthermore, 
\[\gamma(G + \Lb) = (G, 0) + \widetilde{L} = (G, \phi(L)) + \widetilde{L}\] is dense in 
$\mathbb{T}(\Lb)$.
So $\gamma$ is surjective. \qed

\medskip

The proofs of the existence of the torus parametrization shown here are essentially due to Schlottmann
\cite{martin}. There he proves the existence of the map under the condition that a point set $\Lam$ satisfies $\Lam({W}^{\circ}) \subset \Lam \subset \Lam(W)$, where $W$ is compact, and is repetitive (instead of assuming that the window $W$ is the closure of its interior). In fact it is proved in \cite{martin} that the repetitivity of the point set implies that the window is the closure of its interior. However they are not equivalent conditions. 
Here we assume that each of the windows is the closure of its interior and we also place his results into multi-set setting. We just point out here that in the multi-set setting the existence of the torus parametrization can also be proved under the assumption of repetitivity of a  multi-set satisfying $\Lam({W_i}^{\circ}) \subset \Lam_i \subset \Lam(W_i)$, where $W_i$ is compact, without the assumption that each of the windows is the closure of its interior.  

\medskip

Let $\mathbb{X}_g := \{\Gb \in \mathbb{X}(\Lb) : \phi(L) \cap (c_{\Gbt} + \bigcup_{i = 1}^m \partial{W_i}) = \emptyset\}$.
Note that $\mathbb{X}_g \neq \emptyset$ by the Baire category theorem.

\begin{prop} \label{generic}
Let $\Lb$ be a multi-set in $G$. Suppose that $\Lam({W_i}^{\circ}) \subset \Lam_i \subset \Lam(W_i)$,
$W_i$ is compact, and $W_i = \overline{{W_i}^{\circ}} \neq \emptyset$ for all $i \le m$ with respect to 
some irredundant CPS. Then $\gamma |_{\mathbb{X}_g}$ is one-to-one.
\end{prop}

\noindent {\sc Proof.} From Proposition \ref{between-interior-closure} we know that for any 
$\Gb \in \mathbb{X}_g$ there exists $(-s, -c_{\Gbt}) \in G \times H$ such that 
\[-s + \Lam(c_{\Gbt} + {W_i}^{\circ}) \subset \Gam_i \subset -s + \Lam(c_{\Gbt} + W_i) ~~\mbox{for all}~ i \le m.\]
Here $\phi(L) \cap (c_{\Gbt} + \bigcup_{i=1}^m \partial{W_i}) = \emptyset$.
So $\Gam_i = -s + \Lam(c_{\Gbt} + {W_i}^{\circ})$ for all $i \leq m$.
Thus $\gamma |_{\mathbb{X}_g}$ is one-to-one. \qed

\medskip
We define $\mathbb{T}_g := \gamma(\mathbb{X}_g)$.

\begin{prop}
Let $\Lb$ be a multi-set in $G$. Suppose that $\Lam({W_i}^{\circ}) \subset \Lam_i \subset \Lam(W_i)$,
$W_i$ is compact, and $W_i = \overline{{W_i}^{\circ}} \neq \emptyset$ for all $i \le m$, with respect to some irredundant CPS. 
Then $\mathbb{T}_g = \mathbb{T}(\Lb) \backslash ((G \times \bigcup_{i=1}^m \partial{W_i} + \widetilde{L})/ \widetilde{L})$.
\end{prop}
 
\noindent {\sc Proof.} Note that 
\be
(-s, -c_{\Gbt}) + \widetilde{L} \in \mathbb{T}_g 
& \Leftrightarrow &(c_{\Gbt} + \bigcup_{i=1}^m \partial{W_i}) \cap \phi(L) = \emptyset \nonumber \\
& \Leftrightarrow & c_{\Gbt} \notin \phi(L) - \bigcup_{i=1}^m \partial{W_i} \nonumber \\
& \Leftrightarrow & (-s, -c_{\Gbt}) \notin G \times (-\phi(L) + \bigcup_{i=1}^m \partial{W_i}) 
\nonumber \\
& \Leftrightarrow & (-s, -c_{\Gbt}) + \widetilde{L} \notin (G \times \bigcup_{i=1}^m \partial{W_i} + 
\widetilde{L})/\widetilde{L}.
\ee
\qed

\begin{prop} \label{fullmeasure}
Let $\Lb$ be a multi-set in $G$. Suppose that $\Lam({W_i}^{\circ}) \subset \Lam_i \subset \Lam(W_i)$,
$W_i$ is compact, and $W_i = \overline{{W_i}^{\circ}} \neq \emptyset$ for all $i \le m$, with respect to some irredundant CPS. Suppose that $\theta_H (\partial{W_i}) = 0$ for all 
$i \leq m$, where $\theta_H$ is a Haar measure in $H$.
Then $\lambda(\mathbb{T}_g) = 1$, where $\lambda$ is a Haar measure in $\mathbb{T}(\Lb)$.
\end{prop}

\noindent {\sc Proof.} We will show that $\lambda((G \times \bigcup_{i=1}^m \partial{W_i} + \widetilde{L})/ \widetilde{L}) = 0$. By \cite[Theorem 3.3.28]{RS}, we only need to show that 
$\nu(G \times \bigcup_{i=1}^m \partial{W_i} + \widetilde{L}) = 0$, where $\nu$ is 
a Haar measure in $G \times H$. Since $G$ is compactly generated, there exists a sequence of compact sets $\{K_n\}$ such that $G = \bigcup_{i=1}^{\infty} K_n$. Each $\nu(K_n \times \bigcup_{i=1}^m \partial{W_i})= 0$ from the assumption that $\theta_H (\partial W_i) = 0$ for all $i \leq m$. Since $\widetilde{L}$ is countable,
$\nu(G \times \bigcup_{i=1}^m \partial{W_i} + \widetilde{L}) = 0$. Thus the assertion follows. \qed

\medskip

The following proposition is a modification of \cite[Prop.\,5.1]{MS} for multi-sets.

\begin{prop} \label{AequalsT}
Let $\Lb$ be a multi-set in $G$. Suppose that $\Lam({W_i}^{\circ}) \subset \Lam_i \subset \Lam(W_i)$,
$W_i$ is compact, and $W_i = \overline{{W_i}^{\circ}} \neq \emptyset$ for all $i \le m$, with respect to some irredundant CPS. Suppose that $\theta_H (\partial{W_i}) = 0$ for all 
$i \leq m$, where $\theta_H$ is a Haar measure in $H$.
Then $\mathbb{A}(\Lb) \cong \mathbb{T}(\Lb)$.
\end{prop}

\noindent {\sc Proof.} $\mathbb{A}(\Lb)$ is the completion of $G$ under the autocorrelation topology. $\mathbb{T}(\Lb)$ may be considered as the completion of $G$ when it is given the coarsest topology for which the mapping $x \mapsto (x,0) + \widetilde L$ of $G$ into $\T(\Lb)$ is continuous. It will suffice to show that these two topologies on  $G$ are the same. If $x \in G$ is close to $0$ in $\mathbb{T}$-topology, then for small open neighbourhoods $V$ of $0$ in $G$ and $V_1$ of $0$ in $H$ there exists $(t, \phi(t)) \in \widetilde{L}$ such that $x - t \in V$ and $\phi(t) \in V_1$. On the other hand if $x \in G$ is close to $0$ in $\mathbb{A}$-topology, then for a small open neighbourhood $U$ of $0$ in $G$ and some $\epsilon > 0$ there exists $t \in L$
 such that $x - t \in U$ and $d(t+ \Lb, \Lb) < \epsilon$. So we need to show that for $t \in L$, $\phi(t)$ is close to $0$ in $H$ if and only if $d(t + \Lb, \Lb)$ is close to $0$.

For $t \in L$, 
\be
d(t + \Lb, \Lb) &=& \lim_{n \rightarrow \infty} \supreme \frac{\sum_{i=1}^m \sharp(((t + \Lam_i)\; \triangle \; \Lam_i) \cap A_n)}{\theta(A_n)} \nonumber \\
&=& \sum_{i=1}^m \lim_{n \rightarrow \infty} \frac{ \sharp(((t + \Lam_i) \; \triangle \; \Lam_i) \cap A_n)}{\theta(A_n)} \nonumber \\
&=& \sum_{i=1}^m (\theta_{H}(W_i \backslash (\phi(t) + W_i)) + \theta_{H} (W_i \backslash (-\phi(t) + W_i))),  
\ee
since each point set is a regular model set by the assumption (see \cite{RVM1}). 

Note that 
\[\theta_{H} (W_i \backslash (s + W_i)) = \theta_{H}(W_i) -  \mathbf{1}_{W_i} * 
\widetilde{\mathbf{1}_{W_i}}(s)\]
 is uniformly continuous in $s$ (see \cite[Subsec.\,1.1.6]{Rudin}). So if $\phi(t)$ converges to $0$ in $H$, 
then $d(t + \Lb, \Lb)$ converges to $0$ in $\R$. 

On the other hand, suppose that $\{t_n\}$ is a sequence
such that $d(t_n + \Lb, \Lb) \rightarrow 0$ as $n \rightarrow \infty$.
Then for each $i \le m$ 
\[\{\theta_{H} (W_i \backslash (\phi(t_n) + W_i))\}_n \rightarrow 0~ \mbox{as}~ n \rightarrow \infty.
\]
Note that for large enough $n$, $W_i \cap (\phi(t_n) + W_i) \neq \emptyset$ and so 
$\phi(t_n) \in W_i - W_i$ for all $i \leq m$. Since $W_i - W_i$ is compact, $\{\phi(t_n)\}_n$ has a converging subsequence 
$\{\phi(t_{n_k})\}_k$. For any such sequence define ${t_0}^* := \lim_{k \rightarrow \infty} \phi(t_{n_k})$.
Then 
\[\theta_{H} (W_i \backslash ({t_0}^* + W_i)) = 0
\]
and so $\theta_{H} ({W_i}^{\circ} \backslash ({t_0}^* + W_i)) = 0$ for each $i \le m$. 
Thus ${W_i}^{\circ} \subset {t_0}^* + W_i$ and this implies 
$W_i \subset {t_0}^* + W_i$. On the other hand, 
$\lim_{k \rightarrow \infty} - \phi(t_{n_k}) = - {t_0}^*$ and $\theta_{H} ({W_i}^{\circ} \backslash (-{t_0}^* + W_i)) = 0$. 
So $W_i \subset -{t_0}^* + W_i$. Hence 
$W_i \subset {t_0}^* + W_i \subset {t_0}^* - {t_0}^* + W_i$ 
and $W_i = {t_0}^* + W_i$.
This equality is for each 
$i \leq m$. Thus ${t_0}^* = 0$ by irredundancy. So all converging subsequences $\{\phi(t_{n_k})\}_k$ converge to $0$ and
$\{\phi (t_n)\}_n \rightarrow 0 ~\mbox{as}~ n \rightarrow \infty$.

This establishes the equivalence of the two topologies. By \cite[Prop 5, III \S 3.3]{Bour}, there exists an isomorphism of $\mathbb{A}(\Lb)$ onto $\mathbb{T}(\Lb)$. \qed

\medskip
 
Theorem \ref{main}  is a direct consequence of Propositions \ref{continuity}, 
\ref{generic}, \ref{fullmeasure}, and \ref{AequalsT}.    

\section*{Acknowledgements} The authors would like to thank the Banff International Research Station for its hospitality during the initial stages of this work. RVM thanks the Mathematisches Forschung Instit\"ut Oberwolfach for the RIP programme which supported him during part of this research. Finally he would like to thank Michael Baake and Daniel Lenz for their collaboration in much of the mathematics that inspired this article. That work is to be reported in \cite{BLM}.

\end{document}